\documentclass[12pt, a4paper]{article}
	
\usepackage{amsmath, amsthm, amscd, amsfonts, amssymb, graphicx, color}
\usepackage[bookmarksnumbered, plainpages]{hyperref}
	
\begin{document}
 \centerline{\Large{\bf A Construction of Best Fractal Approximation}}
 \centerline{}
 \centerline{Yong Suk Kang, \textsuperscript{a}Chol Hui Yun, Dong Hyok Kim}
 \centerline{}
 \small \centerline{Faculty of Mathematics, \textbf{Kim Il Sung} University,}
\small \centerline{Pyongyang, D.P.R Korea}
\small \centerline{\textsuperscript{a} Corresponding author. e-mail address: yuncholhui@yahoo.com}
\centerline{}
\centerline{}
\begin{abstract}
In this paper we present a method for constructing the continuous best fractal approximation in the space of bounded functions. 
We construct the finite-dimensional subspace of the space of bounded functions whose base consists of the continuous fractal functions, 
and propose how to find the best approximation of given continuous function by element of the constructed space.
\end{abstract}
{\bf Keywords:} Iterated function system (IFS), Recurrent iterated function system (RIFS), Fractal function, Fractal approximation, 
Fractal interpolation, Best approximation \\
{\bf MSC 2010:}  37C45, 28A80, 41A05
%
%
%
%
\section{Introduction}
Fractal functions have been widely used in approximation theory, signal process, computer graphics, modeling of natural objects. 
Fractal approximation is important because there are a lot of objects which must be modeled by fractals  in the nature. 
Therefore, constructions of fractal interpolation functions and fractal approximation have been studied in many papers \cite {ban1}-\cite {zh}.

Fractal interpolation has been used in a construction of fractal functions by many scientists. In \cite {ban1}, a construction of one variable fractal functions by IFS with a data set on $\mathbb{R}$ was presented and it 
was generalized in \cite {ban3}, \cite {bou1} by using RIFS.

Constructions of bivariate fractal interpolation functions(BFIFs) have been considered in many papers. 
A method of constructing BFIFs by fractal interpolation functions was introduced in \cite {bou2}, \cite {yun3},  and a construction of self-affine fractal interpolation functions with data set on triangular domain where the interpolation points on the boundary data are coplanar was presented in \cite {mass}. 
In \cite {dalla}, fractal interpolation functions  were constructed  in the case where interpolation points on the boundary data sets on a grid of 
rectangular are collinear, which were generalized in \cite {mal}, \cite {my1}, where any data set on the rectangular was used. 
The construction of recurrent fractal interpolation function on $\mathbb{R}^n$ was presented in \cite {bou1} and fractal functions with function scaling factors 
were constructed in \cite {feng}, \cite {my1}, \cite {yun2}. 
In \cite {ka}, super iterated function system (SIFS), which is a generalization of the IFS, was used in the construction of more complicated fractal functions.

In \cite {nav}, \cite {zh} construction and best approximation of functions by the fractal functions have been studied.  
The fractal  functions constructed by the interpolation are always  continuous, see \cite {ban1}-\cite {nav},\cite {yun2},\cite {yun3}. 
But since in best approximation one tries to find an approximation which minimizes  $\|f-Tf \|$  in the $L^2$ space by Collage theorem, 
the found approximation is not always continuous \cite {st}, \cite {zh}.
We present a construction of continuous best approximation in   the $L^2$ space.


\section{The space of fractal functions}
Let us denote a closed interval [ $a, b$ ]  by $ \Omega$,  $\Omega  \times \mathbb{R}$ by  $X$ and let $N \in \mathbb{N}$, $N \geq 2$ . 
Let mappings $u_l:\Omega \rightarrow \Omega$, $l=1,\cdots , N-1$  be contraction  bijections with $ \bigcup_{l=0}^{N-1}u_l(\Omega) = \Omega$  
and function $v_l:X\rightarrow R$, $l=1,\cdots , N-1$ be Lipschtz functions with respect to second variable.

\newtheorem{Thot}{Definition}
\begin{Thot}\label{def1}
\textnormal{\cite{bou1}} A function  $f(x)$ defined by
\begin{equation}
f(x)=\sum_{l=0}^{N-1}\left[v_l\left(u_l^{-1}(x), f(u_l^{-1}(x)) \right) \right]\chi _{u_l(\Omega)} \label{eq1}
\end{equation}
is called a fractal function, where  $\chi _{u_l(\Omega)}$ is an indicator function of the set $u_l(\Omega)$.

Then the graph of function $f(x)$  consists of a finite number of copies of itself.
\end{Thot}


\newtheorem{lem1}{Lemma}
\begin{lem1}\label{lemma1}
Let us denote a linear space of the bounded functions on  $\overline{\Omega}$ by $F$ . 
Then $F$ is a complete normed space with respect to the super norm $\| \cdot \|_\infty$  .
\end{lem1}
Now, let  $a=x_0 < x_1 < \cdots< x_N=b$ be a partition of the interval $[a,b]$. We define  $ u_l $, $ v_l $ as follows:             
\begin{eqnarray}
& u_l(x)=a_lx+b_l,~ u_l(x_0)=x_l, ~u_l(x_N)=x_{l+1}, \\  \label{eq2}
& v_l(x,y)=s_ly+\lambda_l(x), ~l=1, \cdots , N-1, \nonumber
\end{eqnarray}
where $s_l$ , $\lambda_l$, $l=1, \cdots, N-1$ satisfy the conditions $|s_l|<1$, $\lambda_l \in F$ . We denote
\begin{eqnarray}
s:=(s_0,s_1,\cdots , s_{N-1}),~~\lambda :=(\lambda_0 , \lambda_1 ,\cdots \lambda_{N-1})\in \prod_{l=0}^{N-1}F  \label{eq3}
\end{eqnarray}
We define an operator $B:F \rightarrow F$  by
\begin{equation}
(Bf)(x)=\sum_{l=0}^{N-1}\left[s_l \cdot f\left(u_l^{-1}(x)\right)+\lambda_l \left(u_l^{-1}(x)\right)\right]\chi _{u_l(\Omega)},  ~~x\in [a,b)  \label{eq4}
\end{equation}
\begin{equation*}
(Bf)(b)=\lim_{x\rightarrow b-0}(Bf)(x)
\end{equation*}
Then the operator  $B$ is a contraction one  on  $F$ and has the unique fixed point $f^{*}(\in F)$, which  satisfies \eqref{eq1}, 
and thus it is the fractal function.
 
In latter of  this paper, we fix  $s$ in \eqref{eq3}. Then the operator $B$ and the fractal function  $f$ satisfying \eqref{eq1} are dependent only on $\lambda$ . 
So we denote them by $B_\lambda$, $f_\lambda$ respectively.

   
\newtheorem{lem2}{Lemma}
\begin{lem1}\label{lemma2}
If $f^{(1)}$, $f^{(2)}$ are the fractal functions corresponding to $\lambda^{(1)}$, $\lambda^{(2)}(\in \prod_{l=0}^{N-1}F )$, 
then for any $\alpha_1$, $\alpha_2 (\in \mathbb{R})$, $\alpha_1f^{(1)}+\alpha_2f^{(2)}$ is the fractal function corresponding to   $\alpha_1\lambda^{(1)}+\alpha_2\lambda^{(2)}$.
\end{lem1}
\begin{proof}
As  $s$ is fixed above, given  $\lambda^{(i)}=(\lambda_0^{(i)}, \lambda_1^{(i)},\cdots , \lambda_{n-1}^{(i)})$, $i=1,2$, 
the operators $B_{\lambda^{(i)}}$, $i=1,2$ are defined by \eqref{eq4}. $B_{\lambda^{(i)}}$, $i=1,2$ are contraction operators on $F$  
and have the unique fixed points $f^{*}_{\lambda^{(i)}}$, $i=1,2$, which are just the fractal functions $f^{(i)}$, $i=1,2$. 
Since
\begin{eqnarray*}
f^{(i)}(x)=\sum_{l=0}^{N-1}\left[s_l \cdot f^{(i)}\left(u_l^{-1}(x)\right)+\lambda^{(i)}_l\left(u_l^{-1}(x)\right) \right]\chi _{u_l(\Omega)}, 
\end{eqnarray*}
we have
\begin{eqnarray*}
\alpha_i f^{(i)}(x)=\sum_{l=0}^{N-1}\left[s_l \cdot  \alpha_i f^{(i)}\left(u_l^{-1}(x)\right)+\alpha_i \lambda^{(i)}_l\left(u_l^{-1}(x)\right)\right]\chi _{u_l(\Omega)}, 
\end{eqnarray*}
thus,
\begin{eqnarray}
\left( \alpha_1 f^{(1)}  +  \alpha_2 f^{(2)} \right)(x) & = & \sum_{l=0}^{N-1}\Big[s_l \cdot  \left(\alpha_1 f^{(1)} +\alpha_2 f^{(2)}\right) \left(u_l^{-1}(x)\right) \nonumber \\
& + & \left( \alpha_1  \lambda^{(1)}_l+\alpha_2  \lambda^{(2)}_l \right) \left(u_l^{-1}(x)\right) \Big] \chi_{u_l(\Omega)}. \label{eq5}
\end{eqnarray}
We define an operator  $B_{\alpha_1  \lambda^{(1)}+\alpha_2  \lambda^{(2)}}:F\rightarrow F$  by
\begin{eqnarray*}
\left(B_{\alpha_1  \lambda^{(1)}+\alpha_2  \lambda^{(2)}}f \right)(x)=\sum_{l=0}^{N-1} \left[s_l\cdot  f + \left(\alpha_1  \lambda^{(1)}_l+\alpha_2  \lambda^{(2)}_l \right) \left(u_l^{-1}(x) \right) \right]\chi_{u_l(\Omega)}.
\end{eqnarray*}
Then the operator is contractive and has a unique fixed point, which is $\alpha_1 f^{(1)}+\alpha_2 f^{(2)} $ by \eqref{eq5}.
That is,  $\alpha_1 f^{(1)}+\alpha_2 f^{(2)} $ is a fractal function corresponding to $\alpha_1  \lambda^{(1)}+\alpha_2  \lambda^{(2)}$. 
\end{proof}

Given $\lambda \in \prod_{l=0}^{N-1}F $ , the fractal function  that is the fixed point of $B_{\lambda}$ is  corresponded  uniquely to $\lambda$ .
Thus,  we denote the fixed point by $f_{\lambda}$  and the set of such  $f_{\lambda}$ by $B$ . 
From Lemma~\ref{lemma2},  $B$ is the linear subspace in $F$.

   
\newtheorem{lem3}{Lemma}
\begin{lem1}\label{lemma3}
 $F$ and  $B$  are linearly isomorphic.
\end{lem1}


\newtheorem{theo1}{Theorem}
\begin{theo1}\label{theorem1}
 Let $f_{\lambda^{(i)}}$, $i=1, \ldots ,n$ ,   be fractal function corresponding to  $\lambda^{(i)}$.  $f_{\lambda^{(1)}}$, 
$f_{\lambda^{(2)}}$,$\cdots$,$f_{\lambda^{(n)}}$ are linearly independent if and only if $\lambda^{(1)}$, $\lambda^{(2)}$, 
$\cdots$, $\lambda^{(n)}$ are linearly independent.
\end{theo1}
\begin{proof}
Let $f_{\lambda^{(1)}}$, $f_{\lambda^{(2)}}$,$\cdots$,$f_{\lambda^{(n)}}$  be linearly independent. 
It is enough to prove that $c_1\lambda^{(1)}$ + $c_2\lambda^{(2)}$ + $\cdots$ + $c_n\lambda^{(n)} = 0$ $\Leftrightarrow $ $c_1=0$, 
$c_2=0$,$\cdots$,$c_n=0$.
Lemma~\ref{lemma2} shows that the fractal function corresponding to $c_1\lambda^{(1)}$ + $c_2\lambda^{(2)}$ + $\cdots$ + $c_n\lambda^{(n)}$ is $c_1f_{\lambda^{(1)}}$ + $c_2f_{\lambda^{(2)}}$ + $\cdots$ + $c_nf_{\lambda^{(n)}}$. 
Since the fractal function  $f(x)$ corresponding to 0 is the function with $f(x)\equiv 0$, if $c_1\lambda^{(1)}$ + 
$c_2\lambda^{(2)}$ + $\cdots$ + $c_n\lambda^{(n)}=0$, then we have $c_1f_{\lambda^{(1)}}$ + $c_2f_{\lambda^{(2)}}$ + 
$\cdots$ + $c_nf_{\lambda^{(n)}}=0$. 
Therefore,  by hypothesis of the theorem $c_1f_{\lambda^{(1)}}$ + $c_2f_{\lambda^{(2)}}$ + $\cdots$ + $c_nf_{\lambda^{(n)}}=0$ 
if and only if $c_1=0$, $c_2=0$,$\cdots$, $c_n=0$ . 
This means that $\lambda^{(1)}$, $\lambda^{(2)}$, $\cdots$, $\lambda^{(n)}$ are linearly independent.

Let  $\lambda^{(1)}$, $\lambda^{(2)}$, $\cdots$, $\lambda^{(n)}$   be linearly independent. 
We assume that $c_1f_{\lambda^{(1)}}$ + $c_2f_{\lambda^{(2)}}$ + $\cdots$ + $c_nf_{\lambda^{(n)}}=0$. 
By Lemma~\ref{lemma2}, we get 
 \begin{eqnarray*}
&& \left(c_1 f_{\lambda ^{(1)}}+ c_2 f_{\lambda ^{(2)}}+ \cdots + c_n f_{\lambda ^{(n)}} \right)(x)=\sum_{l=0}^{N-1} \Big[ s_l \cdot \Big(c_1 f_{\lambda ^{(1)}} + c_2 f_{\lambda ^{(2)}}+\cdots\\
&& \qquad + c_n f_{ \lambda ^{(n)}}\Big) \left(u_l^{-1}(x) \right)+ \left(c_1 \lambda _l^{(1)}+c_2 \lambda _l^{(2)}+ \cdots +c_n \lambda _l^{(n)}\right)\left(u_l^{-1}(x) \right) \Big] \chi_{u_l( \Omega )}.
\end{eqnarray*}
Thus, we have
\begin{eqnarray*}
\sum_{l=0}^{N-1}\left[ \left(c_1\lambda ^ {(1)}_l+ c_2\lambda ^ {(2)}_l+ \cdots + c_n \lambda ^{(n)}_l \right) \left(u_l^{-1}(x) \right) \right] \chi _{u_l( \Omega )}=0, 
\end{eqnarray*}
i.e.  for $ x \in u_l( \Omega ) $ ,
\begin{eqnarray*}
\left(c_1\lambda ^ {(1)}_l+ c_2\lambda ^ {(2)}_l+ \cdots + c_n \lambda ^{(n)}_l)(u_l^{-1}(x)\right)=0. 
\end{eqnarray*}
Therefore, we get $ \left(c_1\lambda ^ {(1)}_l+ c_2\lambda ^ {(2)}_l+ \cdots + c_n \lambda ^{(n)}_l \right)(t)=0 $, $ t \in \Omega $, $ l=0, \cdots ,N-1 $, 
which means that $ c_1\lambda ^ {(1)}+ c_2\lambda ^ {(2)}+ \cdots + c_n \lambda ^{(n)} = 0 $.  
Since $ \lambda^{(1)} $, $ \lambda^{(2)} $, $ \cdots $, $ \lambda^{(n)} $ are linearly independent, $ c_1\lambda ^ {(1)}+ c_2\lambda ^ {(2)}+\cdots + c_n \lambda ^{(n)} = 0 $  if and only if $ c_1=0 $, $ c_2=0 $,$ \cdots $,$ c_n=0 $ . 
Thus,  $ f_{\lambda^{(1)}} $, $ f_{\lambda^{(2)}} $, $ \cdots$, $f_{\lambda^{(n)}} $   are linearly independent.
\end{proof} 

 We can see that the space of fractal functions $B$  has infinite dimension.
If  $ f_{\lambda} \in B$ is a continuous function, then we get
\begin{equation}
\lambda _l(b)-\lambda _{l+1}(a) = s_{l+1}\cdot \frac{\lambda _0(a)}{1-s_0}-s_l\cdot \frac{\lambda_{N-1}(b)}{1-s_{N-1}}, ~ l=0, \cdots , N-1. \label{eq6}
\end{equation}  
In fact, at a point  $x=x_{l+1}, l=0, \cdots, N-1$,
\begin{eqnarray*}
& & f(x_{l+1}-) = s_l \cdot f \left( u_l^{-1}(x_{l+1}) \right) +\lambda_l \left( u_l^{-1}(x_{l+1}) \right) = s_lf(b)+\lambda_l(b) \\
& & f(x_{l+1}+) = s_{l+1} \cdot f \left( u_{l+1}^{-1}(x_{l+1}) \right) +\lambda_{l+1} \left( u_{l+1}^{-1}(x_{l+1}) \right) = s_{l+1}f(a)+\lambda_{l+1}(a),
\end{eqnarray*}
where $f(x_{l+1}+), f(x_{l+1}-)$ denote the right limit and the left limit, respectively.   
By \eqref{eq1}, $f(a)=s_0 \cdot f \left( u_0^{-1}(a) \right) + \lambda_0 \left( u_0^{-1}(a) \right) = s_0f(a)+\lambda_0(a)$. 
Hence $f(a) = \frac{\lambda_0(a)}{1-s_0}$. Similarly, we have $f(b) = \frac{\lambda_0(a)}{1-s_0}$, which gives \eqref{eq6}. 
In the case where $\lambda_l(x), l=0, \cdots, N-1$ are polynomials of degree $d$, we denote the fractal function space corresponding by $B_d$. 
Let us denote a set of all the polynomials of degree less than or equal to $d$ by $P_d$.  
Similarly to Lemma~\ref{lemma3} we can prove that $\prod_{l=0}^{N-1}P_d$ and $B_d$ are linearly isomorphic. 


\newtheorem{theo2}{Theorem }
\begin{theo1}\label{Theorem 2}
Let $\phi_0(x), \phi_1(x), \cdots, \phi_N(x) \in B_1$ be fractal interpolation functions satisfying $\phi_i(x_j)=\delta_{ij}, j=0, 1, \cdots, N$. 
Then $\phi_0, \phi_1, \cdots, \phi_N$ are linearly independent.
\end{theo1}
\begin{proof}
The construction of  fractal interpolation functions gives the result.  
\end{proof}

Let us denote $B'_1:= span \left\{  \phi_0, \phi_1, \cdots, \phi_N \right\}$. Then $B'_1$ is a subspace of $B_1$ of dimension $N+1$ and 
the elements of $B'_1$ are continuous fractal functions.


\section{Construction of best fractal approximation}
Given the continuous function $f$, we consider how to find best approximation element of $f$ by element of $B'_1$ in the sense of $L_2$, 
which is one that find $\left( \alpha_0^*, \alpha_1^*, \cdots, \alpha_N^* \right)$ satisfying the following condition:
\begin{equation}
\underset{(\alpha_0, \cdots, \alpha_N)}{\textnormal{min}} \left\| f-\sum_{k=0}^N \alpha_k \phi_k \right\|_{L_2} = \left\| f-\sum_{k=0}^N \alpha_k^* \phi_k \right\|_{L_2} \label{eq7}.
\end{equation}  


\newtheorem{theo3}{Theorem }
\begin{theo1}\label{Theorem 3}
There exists a unique best fractal approximation satisfying \eqref{eq7}. 
\end{theo1}
\begin{proof}
Since $B'_1 \subset L_2, f \in L_2$ and there is a unique best approximation by elements of space of finite dimension in strictly 
normed space, we get the result.
\end{proof}

To solve \eqref{eq7}, it is enough to solve the normal equation. But since $\phi_0, \phi_1, \cdots, \phi_N$ are fixed points of some contraction mappings, 
it is impossible to get  the exact solutions of \eqref{eq7}. Furthermore, even though we calculate $\phi_0, \phi_1, \cdots, \phi_N$ approximately, 
it is impossible to get their exact representations.  Therefore, we present how to find the best approximation by Lemma ~\ref{lemma2}. 
     
It is obvious that we can define an operator $T:B \to B$ by
\begin{equation}
(Tg)(x) = \sum_{l=0}^{N-1} \left[ s_l \cdot g\left( u_l^{-1}(x) \right) + \sum_{k=0}^{N-1} \alpha_k \lambda_l^{(k)}\left( u_l^{-1}(x) \right) \right] \chi_{u_l}(\Omega) \label{eq8}
\end{equation} 
and by Lemma~\ref{lemma2}  $\sum_{k=0}^N \alpha_k \phi_k$ is the fixed point of  the operator $T$ on $B'_1$, where $\lambda^{(k)}=\left( \lambda_0^{(k)}, \lambda_1^{(k)}, \cdots, \lambda_{N-1}^{(k)} \right), k=0, \cdots, N-1$ are functions corresponding to the function $\phi_k$.

Let us denote the fixed point of the operator $T$ by $f_T$ and a contractive constant of $T$ by $c$. 
Then Collage Theorem shows that if $\|f-Tf\|<\varepsilon$ for $\varepsilon>0$, then $\|f-f_T\|<\frac{\varepsilon}{1-c}$. 
On the basis of the theorem, we find an operator, which is denoted by $T^*$, such that $\|f-Tf\|$ is minimum.  
\begin{eqnarray*}
\|f-Tf\|_{L_2}^2 & = & \sum_{l=0}^{N-1} \int_{x_l}^{x_{l+1}} \left[ s_l \cdot f \left( u_l^{-1}(x) \right) + \sum_{k=0}^N \alpha_k \lambda_l^{(k)}
\left( u_l^{-1}(x) \right)-f(x) \right]^2 dx\\
	& = & \sum_{l=0}^{N-1} a_l \int_a^b \left[s_l \cdot f(x) + \sum_{k=0}^N \alpha_k \lambda_l^{(k)}(x)-f(u_l(x)) \right]^2 dx \\
	& = & \sum_{l=0}^{N-1} a_l \left( s_lf+\sum_{k=0}^N \alpha_k \lambda_l^{(k)}-f \circ u_l, ~~ s_lf+\sum_{k=0}^N \alpha_k \lambda_l^{(k)}-f \circ u_l \right),
\end{eqnarray*}
which is denoted by $\phi(\alpha_0, \alpha_1, \cdots, \alpha_N)$. To minimize $\phi(\alpha_0, \alpha_1, \cdots, \alpha_N)$, we have to solve the 
equations $\frac{\partial\phi}{\partial\alpha_k}=0, ~~ k=0, \cdots, N$, i.e.

\begin{eqnarray}
&& \sum_{l=0}^{N-1} a_l \left( 2\left( s_lf, \lambda_l^{(k)} \right) + 2\sum_{j=0}^N \alpha_k \left( \lambda_l^{(k)}, \lambda_l^{(j)}\right) - 
2\left(\lambda_l^{(k)}, f \circ u_l\right) \right)  \nonumber \\
&& \qquad \qquad \times \sum_{j=0}^{N-1}\alpha_j \sum_{l=0}^{N-1}a_l \left( \lambda_l^{(k)}, \lambda_l^{(j)} \right) \nonumber \\
&& \qquad = \sum_{l=0}^{N-1} a_l \left(-\left(s_lf, \lambda_l^{(k)}\right) + \left( \lambda_l^{(k)}, f \circ u_l \right) \right), k=0, 1, \cdots, N \label{eq9}
\end{eqnarray} 

Let us denote $\alpha = (\alpha_0, \alpha_1, \cdots, \alpha_N)^T$ and represent \eqref{eq9} by the form of $A\alpha = \beta$, where $A$ is as follows:
\begin{equation*}
A=
\begin {bmatrix}
\sum_{l=0}^{N-1} a_l \left( \lambda_l^{(0)}, \lambda_l^{(0)} \right) & \sum_{l=0}^{N-1} a_l \left( \lambda_l^{(0)}, \lambda_l^{(1)} \right) & 
\cdots & \sum_{l=0}^{N-1} a_l \left( \lambda_l^{(0)}, \lambda_l^{(N)} \right) \\
\sum_{l=0}^{N-1} a_l \left( \lambda_l^{(1)}, \lambda_l^{(0)} \right) & \sum_{l=0}^{N-1} a_l \left( \lambda_l^{(1)}, \lambda_l^{(1)} \right) & 
\cdots & \sum_{l=0}^{N-1} a_l \left( \lambda_l^{(1)}, \lambda_l^{(N)} \right) \\
\cdots & \cdots & \cdots & \cdots \\
\sum_{l=0}^{N-1} a_l \left( \lambda_l^{(N)}, \lambda_l^{(0)} \right) & \sum_{l=0}^{N-1} a_l \left( \lambda_l^{(N)}, \lambda_l^{(1)} \right) & 
\cdots & \sum_{l=0}^{N-1} a_l \left( \lambda_l^{(N)}, \lambda_l^{(N)} \right)
\end {bmatrix}.
\end{equation*}
 
For example, in the case of $a=0, b=1$, the matrix $A$ is calculated as follows:
\begin{eqnarray*}
&& \lambda_l^{(0)}(x) = \left\{
\begin {array}{ll}
(s_0-1)(x-1), & l=0  \\
s_l(x-1), & otherwise
\end {array} \right. \\
&& \lambda_l^{(i)}(x) = \left\{
\begin {array}{ll}
x, & l=i-1  \\
-x+1, & l=i \\
0, & otherwise
\end {array} \right.
\end{eqnarray*}
\begin{eqnarray*}
&& \lambda_l^{(N)}(x) = \left\{
\begin {array}{ll}
(1-s_{N-1})x, & l=N-1  \\
-s_lx, & otherwise
\end {array} \right. \\
&& a_l=\frac{1}{N}
\end{eqnarray*}

%
%


\begin{figure}[tbhp]
\centering
 \includegraphics[width=122mm]{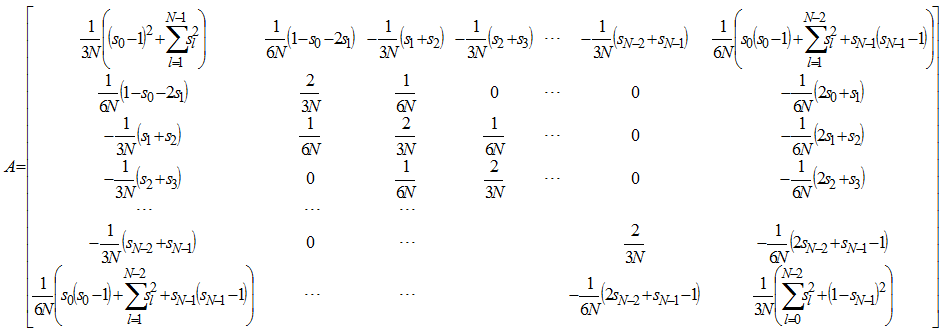}. \\ 
\end{figure}

\vspace{0.5cm}


\noindent Especially,  in the case of $s=0$ we have
\begin{equation*}
A=
\begin {bmatrix}
\frac{1}{3N} & \frac{1}{6N} & 0 & \cdots & \cdots & 0 \\
\frac{1}{6N} & \frac{2}{3N} & \frac{1}{6N} & 0 & \cdots & 0 \\
\cdots & \cdots & \cdots & \cdots & \cdots & \cdots \\
0 & 0 & 0 & \cdots & \frac{1}{6N} & \frac{1}{3N}
\end {bmatrix}.
\end{equation*}
It is clear that the matrix $A$ is non-degenerate, which gives a uniqueness and existence of our best approximation solution.

Finding $\alpha^*$ which is the solution of \eqref{eq9} gives the expected fractal approximation 
\begin{equation}
\sum_{k=0}^N \alpha_k^* \phi_k(x), \label{eq10}
\end{equation}           
which is just the fixed point of  the operator $T^*$. Therefore it is enough that   when finding \eqref{eq10},  we calculate  the fixed point of $T^*$. 


\section{Conclusion}
In this paper we firstly provided a theoretical basis for construction of a finite-dimensional space which consists of continuous fractal functions.  Next we presented a method how to calculate approximately the best approximation by the element of the constructed space. 
The base of the constructed finite-dimensional space consists of the fractal interpolation functions. 
When calculating the best approximation, we use only coefficients of the contraction operator whose fixed points are those fractal interpolation functions instead of interpolation functions themselves. 
Therefore, the presented method has advantage that the calculation of the best approximation is comfortable and the obtained best approximation is continuous.\\

\textbf{Acknowledgement}. We would like to thank editors and referees for their helpful comments and suggestions.

\end{document}